\documentclass{article}
\pdfoutput=1

\usepackage{amstext,amsfonts,amsmath}
\usepackage{graphicx}
\usepackage{doi}
\newcommand{\ud}{\mathrm{d}}
\newcommand{\du}{\, \mathrm{d}}

\newcommand{\R}{\mathbb{R}}

\newcommand{\tol}{\mathrm{tol}}
\newcommand{\tfin}{t_{\textrm{final}}}

\usepackage{subfigure}

\title{Error Control for Exponential Integration of the Master Equation}

\author{Katharina Kormann\thanks{Technische Universitat M\"unchen, Boltzmannstr. 3, D-85747 Garching, Germany ({\tt katharina.kormann@tum.de}).}
        \and Shev MacNamara\thanks{The School of Mathematics and Statistics, University of New South Wales, UNSW, Sydney, Australia {\tt (s.macnamara@unsw.edu.au)}}}

\begin{document}

\maketitle

\begin{abstract}
Error estimates for the numerical solution of the master equation are presented.
Estimates are based on adjoint methods.
We find that a good estimate can often be computed without spending computational effort on a dual problem.
Estimates are applicable to both settings with time-independent, and time-dependent propensity functions.
The Finite State Projection algorithm reduces the dimensionality of the problem and time propagation is based on an Arnoldi exponential integrator, which in the time-dependent setting is combined with a Magnus method.
Local error estimates are devised for the truncation of both the Magnus expansion and the Krylov subspace in the Arnoldi algorithm.
An issue with existing methods is that error estimates for truncation of the state space depend on measuring a loss of probability mass in a way that is not usually compatible with the approximation of the exponential.
We suggest an alternative error estimate that is compatible with a Krylov approximation to the matrix exponential.
Finally, we apply the new error estimates to develop an adaptive simulation algorithm.
Numerical examples demonstrate the benefits of the approach.
\end{abstract}

\textbf{Keywords}: Master equation,  Arnoldi method, adjoint method, chemical kinetics

\textbf{AMS}: 65F10, 65F60, 65L20, 65L70

\pagestyle{myheadings}
\thispagestyle{plain}
\markboth{K. KORMANN AND S. MACNAMARA}{ERROR CONTROL FOR THE MASTER EQUATION}

\section{Introduction}

Markov processes offer a common mathematical framework with which to model applications as diverse as economics, chemical kinetics, single-molecule experiments, ecology, statistical mechanics, telecommunications, complex networks and more \cite{KouKou03,SmiOnn09,ArkRos98,ElfLi07,HilPau11, BlyMcK07,FroMel94,PetzoldCloud,CaoLi04}.
Estimating the probability of being in a certain state is fundamental.
Monte Carlo simulation of trajectories  is one way to compute such estimates   \cite{CoxMil65}.
For example, in the context of chemical kinetics, the Gillespie Stochastic Simulation Algorithm and its variations are popular \cite{Gil77,Kur80}.

Another approach is to directly solve a master equation (sometimes known as a forward Kolmogorov equation), which governs the evolution of the probability distribution associated with the discrete-state, continuous-time Markov process ~\cite{vanKampen07}.
The master equation may be written as a linear system of ordinary differential equations
\begin{equation}\label{eq:master}
	\frac{\du}{\du t} p = A(t) p, \qquad p(0) = p_0.
\end{equation}
Here $p \in \mathbb{R}^N$ is a probability vector: $p \geq 0$ component-wise, and $||p||_{\ell_1}=1$.
The $i$th component is the probability of state $i$.
The off-diagonals of the propensity matrix $A \in \mathbb{R}^{N \times N}$ are nonnegative.
The diagonals are such that each column of $A$ sums to zero.
The off-diagonals have the interpretation that  $a_{ij} \ud t$ is the probability that flows from state $j$ to state $i$ in the next infinitesimal time interval $\ud t$ \cite{Nor97,Ste94}.
A simple example is the graph Laplacian on a line of nodes \cite{StrMac14}.
This article is concerned with the numerical solution of equation \eqref{eq:master} when the matrix $A$ has this special Markov structure.

When $A$ is constant, finite and bounded, the solution to \eqref{eq:master} is the matrix exponential,
\begin{equation}\label{eq:Master_solution}
       p(t) = \exp(t A ) p(0).
\end{equation}
In some applications the number of states, $N$, and the size of the associated matrix, may be very large.
This presents a computational challenge.
Monte Carlo approaches are popular partly because of this challenge.
Nevertheless, methods for the direct solution of the master equation have also been proposed, together with various approaches to the assessment of accuracy  \cite{DeuWul89, HegBur06, DeuHui07,JahWil08,HerWol09,HelLot07,Munsky06,MacBur08, MacNamara10,WolGoe10}.

One approach to this challenge is the Finite State Projection (FSP) algorithm \cite{Munsky06}, which is related to the method of finite sections \cite{Lindner06}.
A key idea is that, starting with an initial distribution of finite support, most of the probability stays localised nearby when propagating over a small time interval.
Instead of working on the full state space, an FSP algorithm focuses computational effort on a subset of the full state space, which is smaller and therefore more efficient, but which still captures most of the probability, and so is still very accurate.
With the assumption of exact matrix exponentials, and the law of conservation of probability, very attractive \textit{a posteriori} error estimates have been derived \cite{Munsky06}.
However, usually an approximation to the matrix exponential in  (\ref{eq:Master_solution})  is employed,  such as an exponential integrator \cite{Moler03, MacBur08, MacNamara10,Hochbruck99,SchTre07,Sid98,KasTre05,MohHig11}.
These approximations do not exactly conserve probability so obtaining an estimate of the overall error  when combining two or more such  approximations is an issue.
This issue serves as motivation for the present article, which shows how adjoint methods offer a mathematical framework that is more accommodating of combinations of approximations and that can be used for deriving overall error estimates.

When efficient estimates of the \textit{a posteriori} error are available, adaptive step sizes can be used to compute the solution in an efficient way, and in a way that controls the global error.
Adjoint methods find wide applications in engineering, and for ordinary differential equations \textit{a posteriori} error estimates based on the adjoint method have been devised \cite{Cao04,CaoLiPetzoldSerban03,BuiWillcoxGhattas08,AscherPetzoldBook}.
In the context of the time-dependent Schr{\"o}dinger equation, it has been shown that efficient error estimates can come from combining adjoint methods with conservation properties \cite{Kormann11}.
Here, we show how to adapt adjoint approaches such as appear in \cite{Cao04,CaoLiPetzoldSerban03,AscherPetzoldBook} to the master equation.

The outline of the article is as follows.
In the next section, we devise \textit{a posteriori} error estimates for the master equation.
The estimates are based on a duality argument.
Three different levels of approximation, with successively less dependence on the dual problem, are derived.
These are general purpose estimates that may be combined with various numerical methods.
Sec.~\ref{sec:statespace} discusses truncation of the state space.
We take special care of the fact that the matrix exponential cannot be considered to be computed exactly for numerical propagation.
Next, we introduce a Magnus--Arnoldi propagator.
Sec.~\ref{sec:full_estimate} discusses how to compute a discretized error estimate for the chosen numerical integrator.
Sec.~\ref{sec:adaptive} details one way to adaptively solve the master equation based on the \textit{a posteriori} estimates previously introduced.
We discuss the influence of the dual problem on the estimates, and compare different levels of approximation in Sec.~\ref{sec:dual}.
Sec.~\ref{sec:numerical} presents numerical results, with examples from chemistry and  immunology.

\section{\textit{A posteriori} error estimates}\label{sec:apost}

In this section, we devise \textit{a posteriori} error estimates based on duality (cf.~\cite{Cao04,Moon03}).
These are generic estimates that may be combined with various numerical methods.
Sec.~\ref{sec:full_estimate} shows how to tailor these estimates to a specific numerical method based on Magnus--Arnoldi propagation and state space reduction.

The global error is estimated by a weighted sum of the local errors at each step.
The estimate is of the form \;\; error  $\leq$ (condition number) $\times$ (perturbation to data). \;
This form highlights an analogy with the condition number, $\kappa_A$, of  a matrix for a linear system, $Ax=b$, in which the relative error, $\rho_x$, in the solution $x$ is bounded by $\rho_x \leq \kappa_A \rho_b$, where $\rho_b$ is the relative error in (or perturbation to) the data $b$.
In this analogy, the solution of a dual problem, which may be thought of as a sensitivity and which we denote by $q$ below, plays a role similar to the condition number of a matrix, and the local numerical errors, denoted $r$ below, play a role similar to perturbations to the data.

Consider the master equation (\ref{eq:master}).
Let $\widetilde p$ denote a continuation of the discrete solution obtained from a numerical scheme that solves the perturbed equation
\begin{equation}
	\frac{\du}{\du t} \widetilde p = A(t) \widetilde p + r(t), \qquad \widetilde p(0) = p_0 + R, \nonumber
\end{equation}
where $r(t)$ is the residual due to numerical error and $R$ a perturbation in the initial value. As in \cite{Cao04}, we assume that the perturbations are bounded by a small constant.
The error is
\[
E=p(\tfin)-\widetilde p(\tfin)
\]
at time $\tfin$, where $p$ is the exact solution to \eqref{eq:master}. We want to estimate some functional of the error.
Let $z$ be a normalized vector defining this functional.
In order to estimate the functional $z^TE$, we consider the adjoint master equation
\begin{equation}\label{eq:adjoint}
	\frac{\du}{\du t} q = - A^T(t) q, \qquad q(\tfin) = z.
\end{equation}
Note that the eigenvalues of the propensity matrix $A(t)$ all have non-positive real part.
Therefore the forward problem \eqref{eq:master} is stable.
Then, also the adjoint equation is stable, as shown in \cite[Thm.~4.1]{CaoLiPetzoldSerban03}.
The fundamental theorem of calculus relates our error functional to the perturbation $r(t)$,
\begin{equation}\begin{aligned}\label{eq:est}
z^T E &= \int_0^{\tfin} (\frac{\du}{\du t} q)^T E + q^T(\frac{\du}{\du t} E) \du t + q(0)^T E(0) \\&= \int_0^{\tfin} (-A(t)^T q)^T E + q^T(A(t) p(t)-\tilde p(t) - r(t)) \du t + q(0)^T R \\&= -\int_0^{\tfin} q(t)^T r(t) \du t + q(0)^T R.
\end{aligned}\end{equation}
Apply H{\"o}lder's inequality to estimate (\ref{eq:est}) by
\begin{equation}\label{eq:est1}
	|z^T E |   \leq \int_0^{\tfin}\|q(t)\|_{\infty} \|r(t)\|_{\ell_1} \du t + \|q(0)\|_{\infty}\|R\|_{\ell_1}.
\end{equation}
Many applications give rise to self-adjoint problems and $\ell_2$-norms are a natural choice, but the $\ell_1$-norm and $\ell_{\infty}$-norm are more suitable in our setting of the Master equation.

In order to approximate (\ref{eq:est}), consider the conservation properties of the adjoint equation (\ref{eq:adjoint}), which can be deduced from the following properties of the master equation (cf.~\cite{vanKampen07}):
\begin{eqnarray}
	\sum_{k} A_{kl} &=& 0  \quad \textrm{for each } l, \label{eq:zero:column:sum}\\
	 A_{kl} &\geq& 0 \quad \textrm{for } k\neq l. \label{eq:nonnegative:off:diagonal}
\end{eqnarray}
The first property gives
\begin{equation}\label{eq:adjoint2}
	\frac{\du}{\du t}q_k = -\sum_{l} (q_l-q_k) A_{lk}.
\end{equation}
Let $k_*$ be the index with $q_{k_*}(t) = \max_l q_l(t)$ at a certain point in time.
The off-diagonals of $A$ are nonnegative so from (\ref{eq:adjoint2}) we deduce $\frac{\du}{\du t}q_{k_*} >0$.
Hence, the $\ell_{\infty}$ norm is decreasing when going backwards in time, i.e.,
\begin{equation}\label{eq:cons}
	\max_{t \in [0,\tfin]} \|q(t)\|_{\infty} = \|q(\tfin)\|_{\infty} = \|z\|_{\infty} .
\end{equation}
This further estimates $|z^TE|$ by
\begin{equation}\label{eq:est2}
	|z^T E |    \leq \|z\|_{\infty} \left(\int_0^{\tfin} \|r(t)\|_{\ell_1} \du t+\|R\|_{\ell_1} \right).
\end{equation}

\paragraph*{An estimate without numerical evaluation of the dual problem}
At this point, we have derived three \textit{a posteriori} error estimates: \eqref{eq:est}, \eqref{eq:est1}, and \eqref{eq:est2}.
We make some observations about their relative merits.
Usually, a disadvantage of adjoint methods is that they require the solution to the dual problem, which is expensive both to compute and to store in memory.
For example, to apply estimate \eqref{eq:est}, we first need to solve the dual problem \eqref{eq:adjoint} for $q(t)$ and save this solution in memory for $(0, \tfin)$, so that we can then evaluate the integral in \eqref{eq:est}.
In contrast, note that estimate (\ref{eq:est2})  does \emph{not} require the solution of a dual problem.
This makes estimate (\ref{eq:est2}) very attractive.
Compared to  estimate (\ref{eq:est}), estimate (\ref{eq:est2}) is cheaper to compute but may overestimate the error.
Whether or not it pays off to solve the dual is problem and parameter dependent.
Below we show that  (\ref{eq:est2}) is particularly beneficial if we want to bound several functionals of the error at the same time.

\paragraph*{Remark}
The result of equation (\ref{eq:cons}) is based on the property $\sum_{k} A_{kl} = 0$ for each $k$.
Truncating the state space, as described above, destroys this property.
Then, estimates (\ref{eq:est1}) and (\ref{eq:est2}) hold all the same.
Instead of (\ref{eq:adjoint2}), we have
\begin{equation}
	\frac{\du}{\du t}q_k = -\sum_{l \in \mathcal{I}_1} (q_l-q_k) A_{lk} + \sum_{l \in \mathcal{I}_2} q_k A_{lk}. \nonumber
\end{equation}
Since $A_{lk} \geq 0$, $l \neq k$, and $q_k \geq 0$, we still have $\frac{\du}{\du t} q_{k^*} > 0$ for the index $k^*$ where $p$ takes its maximal value.
This was the requirement for estimates (\ref{eq:est1}) and (\ref{eq:est2}).

\paragraph*{A component-wise error bound}
Bounding the component-wise error is of interest \cite{Munsky06}.
For this purpose, consider functionals $z$ defined by standard basis vectors $z=e_1, z=e_2, \ldots$, where $e_i$ is the $i$th column of the identity matrix, and observe that all have an $\ell^{\infty}$ norm of $1$.
Hence, by \eqref{eq:est2} the error in each component is bounded by $\int_0^{\tfin} \|r\|_{\ell_1} \du t+\|r(0)\|_{\ell_1}.$
In this way, we have found a component-wise error bound that does not require knowledge of the adjoint equation.
In that case, for a given tolerance $\varepsilon>0$, if we control the error so that
\[
 \int_0^{\tfin} \|r(t)\|_{\ell_1} \du t + \|R\|_{\ell_1}\leq \varepsilon,
\]
then the exact solution $p$ is within $\varepsilon$ of the numerical solution $\widetilde p$:
\begin{equation}\label{eq:simple:error:bound}
	\widetilde p - \varepsilon e \leq p \leq \widetilde p + \varepsilon e,
\end{equation}
where $e = (1,\ldots,1)$.
This estimate may be compared to estimate (2.9) in \cite{Munsky06} or estimate (2.4) in \cite{MacNamara10}, which is reproduced here in \eqref{eq:FSP:est} for convenience.
However, those approaches focus on one primary source of error, rather than offering strategies to accommodate all sources of error. 
The approach presented here augments previous approaches by providing a unifying framework for including errors that arise from different sources.
One common source of error comes from adaptively tracking the support, as we describe next.

\section{Truncation of the state space}\label{sec:statespace}
In many applications of master equations, the number of all possible states is enormous, whilst in practice, at any particular instant in time, only a small subset of these states contain most of the probability.
This observation motivated numerical methods that adaptively track the support of the probability distribution, reducing memory requirements and computational time.

One approach to truncating the state space is as follows \cite{Munsky06}.
Let $A = A_T+A_R$ where $A_T$ belongs to the truncated state space and $A_R$ is the remainder.
Hence, $A_T$ consists of the square block indexed by the included states only, while $A_R$ collects the rows and columns belonging to the dismissed states.
Let $\mathcal{I}_1$ denote the set of states that are included, and let $\mathcal{I}_2$ denote the states that are not included.
Then $A_T$ is the submatrix of $A$ with rows and columns indexed by the states of $\mathcal{I}_1$.
In MATLAB notation, $A_T = A(\mathcal{I}_1,\mathcal{I}_1)$.
Instead of solving with the full matrix $A$, we solve
\begin{equation}
	\frac{\du }{\du t} p^E(t) = A_T(t) p^E(t) = A(t) p^E(t) - A_R(t) p^E(t). \nonumber
\end{equation}
Assuming exact exponentials, it has been shown \cite{Munsky06} that $\|p^E(t)\|_{\ell_1}$ decreases, and the amount by which it decreases from $1$  is a measure of the component-wise error.
Briefly, if $\|p(0)\|_{\ell_1} = 1$ and $\|p^E(t)\|_{\ell_1} = 1-\delta$, then
\begin{equation} \label{eq:FSP:est}
 p^E  \leq p \leq  p^E + \delta e,
\end{equation}
where $p$ is the exact solution.

An issue with this estimate is that we usually do not compute the matrix exponential exactly.
Therefore, the $\ell_1$ norm of the vector $p$ will not be exactly conserved by the numerical propagator (see also the discussion in Sec.~\ref{sec:discussion}) so there is an issue with a state-space-truncation criterion of the form \eqref{eq:FSP:est}.
To address this issue, observe that on a small time interval $[t,t+\Delta t]$, we can estimate the probability mass that exits the states in $\mathcal{I}_1$ by
\begin{equation}\label{eq:res_sst}
	\int_t^{t+\Delta t} A_R(\tau) p^E(\tau) \du \tau.
\end{equation}
Now, $p^E(t)$ is nonzero only  for states included in $\mathcal{I}_1$.
Hence, computing $A_R p^E(t)$ is cheap because we only need to consider states in $\mathcal{I}_2$ that can be reached from $\mathcal{I}_1$ in one step.
Sec.~\ref{sec:CME} compares this criterion \eqref{eq:res_sst} with \eqref{eq:FSP:est}, using an example in which the numerical method conserves the $\ell_1$ norm.
This shows good agreement between \eqref{eq:res_sst} with \eqref{eq:FSP:est}. Hence, \eqref{eq:res_sst} provides a viable alternative criterion to estimate the error due to state space truncation when the estimation criterion \eqref{eq:FSP:est} is unreliable because the numerical propagator does not exactly conserve the $\ell_1$ norm.


\section{A Krylov approximation to the matrix exponential}
Solving the master equation \eqref{eq:master} requires an approximation to the matrix exponential.
In applications where the matrix is large and sparse (a common case in chemical master equations for example), the action of the exponential on a vector may be computed via a Krylov method, as we now describe \cite{Sid98}.

Over one small time step $\Delta t$, consider the approximation of
\begin{equation}
	p(t + \Delta t) = \exp\left( \Omega \right) p(t).
\end{equation}
For master equations \eqref{eq:master} with constant matrix $A$ we have $\Omega = \Delta t \cdot A$ in mind.
The case of time-varying $A(t)$ requires a more complicated choice of $\Omega$, as we describe in the next section.

We now describe the Arnoldi method. For the purpose of this subsection we consider $p$ to be the `correct' function because we want to isolate the error solely due to the Arnoldi method.
The Arnoldi process generates an orthogonal basis $V_{s} = [v_1 = p(t), v_2, \ldots, v_{s}]$ of the $s$th order Krylov subspace $\mathcal{K}_{s} = \mathrm{span} \{ p(t), \Omega p(t), \ldots, \Omega ^{s-1} p(t)\}$ and an upper Hessenberg matrix $H_{s}$ that represents the matrix $\Omega^{[m]}$ in the Krylov subspace,
\begin{equation}\label{eq:Krylov}
	\Omega V_{s} = V_{s} H_{s} + [H_{s+1}]_{s+1,s} v_{s+1}.
\end{equation}
Here, $ [H_{s+1}]_{s+1,s}$ denotes element $(s+1,s)$ of the projection matrix to the $(s+1)$th Krylov subspace and $v_{s+1}$ is the $(s+1)$th vector of the basis $V_{s}$.
We then approximate $p(t+\Delta t)$ by
\begin{equation}
	p(t+\Delta t) \approx p^A(t+\Delta t) = V_{s} \cdot \exp\left(H_{s}\right)\cdot e_1 \cdot \| p(t)\|_{\ell_2},
\end{equation}
where $e_1 = (1,0,\ldots,0) \in \R^{s}$.

Error control and stopping criteria for Krylov subspace methods have been studied by Saad \cite{Saad92} and Hochbruck \& Lubich \cite{Hochbruck99}.
They consider the generalized residual, $\rho_A$, defined as
\begin{equation}\label{eq:def_resL}
	\rho_A = \left(\exp\left( \Omega \right)\right)^{-1} p^A(t+\Delta t)-p(t) .
\end{equation}
From (\ref{eq:Krylov}) they derive the following easy-to-compute expression for the residual:
\begin{equation}\label{eq:arnoldi_residual}
	\rho_{A} = \|p(t)\|_{\ell_2} \Delta t L^{s}_{s+1,s}[\exp\left(L^{s}\right)]_{s,1}\cdot v_{s+1}.
\end{equation}
To see how the residual $\rho_A$ fits into \textit{a posteriori} estimate (\ref{eq:est}), we rearrange Eqn.~(\ref{eq:def_resL}),
\begin{equation}
	p^A(t+\Delta t) = \exp\left(\Omega  \right) \left(p(t)+\rho_A\right),
\end{equation}
i.e. we solve for the exact right-hand-side, but with a perturbed initial value.

\section{A Magnus integrator for a master equation with time-varying rates}
\label{sec:discretization}
When the matrix $A$ is constant, the solution of \eqref{eq:master} is the matrix exponential (\ref{eq:Master_solution}).
When $A=A(t)$ is time-dependent, the evolution operator is more complicated.
For sufficiently small intervals $\Delta t$,
\begin{equation}\label{eq:magnus}
	p(t+\Delta t) = \exp\left(\Omega(t)\right)p(t),
\end{equation}
where the evolution operator is given by the Magnus expansion \cite{Magnus54, Blanes09}, of which the first few terms are:
\begin{equation}
	\Omega (t) = \int_{t}^{t+\Delta t} A(\tau) \du \tau - \frac{1}{2} \int_{t}^{t+\Delta t} \left[\int_{t}^{\tau}A(\sigma) \du \sigma, A(\tau)\right] \du \tau + \ldots. \nonumber
\end{equation}
Here $[A,B]=AB - BA$ is the matrix commutator. All higher order terms in the expansion involve commutators.
In the special case that the matrix commutes with itself, i.e. $[A(t_1), A(t_2)] =0$, those commutators are all zero so the solution to (\ref{eq:master}) is $p(t+\Delta t)=\exp(\Omega(t))p(t)$ where $\Omega (t) = \int_t^{t+\Delta t} A(\tau) \du \tau $.
An example is when $A(t) = c(t) A_c$ for a scalar function $c(t)$ and a constant matrix $A_c$, in which case the integral in $\Omega (t) = A_c (\int_t^{t+\Delta t} c(\tau) \du \tau) $ is scalar.

For a numerical integrator, the expansion is truncated to a certain order in $\Delta t$.
Denoting the truncated Magnus matrix by $\Omega^{[m]}(t)$ with $m$ being the order of the truncation, we propagate based on the formula,
\begin{equation}
	p^M(t+\Delta t) = \exp( \Omega^{[m]}(t)) p^M(t). \nonumber
\end{equation}
Na\"ive computation of successive terms in the Magnus expansion may be expensive.
Fortunately, methods for the efficient computation of the terms in the truncated Magnus expansion up to a certain order have been devised by Iserles and Norsett \cite{Iserles99}.
Second-order and forth-order methods are most common in practice.
Let $\theta_2$ denote the second order truncation and $\theta_4$ the additional term in a fourth order approximation.
We use the optimized formulas presented by Blanes \emph{et al.} \cite{Blanes00}:
\begin{equation}\begin{aligned}
	\theta_{2} = \Delta t B^{(0)}, \qquad \theta_4 = -(\Delta t)^2 [B^{(0)},B^{(1)}], \nonumber
\end{aligned}
\end{equation}
where
\begin{equation}
	B^{(j)}(t_n) = \frac{1}{(\Delta t)^{j+1}}\int_{-\Delta  t/2}^{\Delta t/2} \tau^j A(t+\Delta t/2+\tau) \du \tau, \quad j=0,1,\ldots. \nonumber
\end{equation}

\paragraph*{Error estimates for the truncated Magnus expansion}
Although finer theoretical estimates are known, it is usually more practical to assume that the first omitted term in the expansion dominates the error.
If we propagate based on $\Omega^{[m]}$, the first order error term is
\begin{equation}\label{eq:magnus_residual}
	r_M(\tau) \approx \dot \theta_{m+2} (t) p^{M}(t), \quad t \leq \tau \leq t+\Delta t.
\end{equation}
The residual $r_M$ on interval $[t, t+\Delta t]$ enters the error estimate as
\begin{equation}
	\int_{t}^{t+\Delta t} q(\tau)^T r_M(\tau) \du \tau  \approx \int_{t}^{t+\Delta t} q(\tau)^T \dot \theta_{m+2} (\tau) p^{M}(\tau)   \du \tau  .
\end{equation}
Since we do not know the values of $p(\tau)$ or $q(\tau)$ at any time but $t$ and $t+\Delta t$ (or actually only for $t$ when predicting $\Delta t$), we approximate the integral as
\begin{equation}\begin{aligned}
	\int_{t}^{t+\Delta t} q(\tau)^T \dot \theta_{m+2} (\tau) p^{M}(\tau)   \du \tau  &\approx q(t)^T \int_t^{t+\Delta t} \dot \theta_{m+2}(\tau) \du \tau p^{M}(t) \\
	&=  q(t)^T  \theta_{m+2}(t)  p^{M}(t) .
\end{aligned}\end{equation}
This approximation is sufficiently accurate for our purpose because $\Delta t$  is small and $r_M$ is only approximated to first order.
Finally, we can apply Holder's inequality to the approximated integral for estimates of type (\ref{eq:est1}) or (\ref{eq:est2}).

\paragraph*{Remark} Recall the two properties \eqref{eq:zero:column:sum} and \eqref{eq:nonnegative:off:diagonal} of a matrix $A(t)$ for the master equation.
The commutator $[A(t_1), A(t_2)]$ shares property \eqref{eq:zero:column:sum} with $A(t)$ that columns sum to zero, but unlike $A(t)$, the commutator can have negative off-diagonal entries.

\paragraph*{A special structure of master equations}
In general, the main costs of evaluating approximations to the Magnus expansion, such as the formulas for $\theta_4$, come from the evaluation of the commutators.
However, in the context of master equations there is often a special structure to the time dependence, and this structure can help.
The propensity matrix often has the form
\begin{equation}
A(t) = A_c + \sum_{l=1}^r f_l(t) A_l
\label{eq:special:form}
\end{equation}
where $A_c$ and $A_l$ are constant in time and the $f_l: \R \rightarrow \R$ are scalar functions.
This special case arises when the propensity is of the \textit{separable} form $\alpha(x,t)=\alpha(x) \alpha(t)$. The rate is separable into a product of two terms, one that depends only on the state, and the other that depends only on time.
For example if the bimolecular reaction $A + B \stackrel{c(t)}{\rightarrow} C$ has time varying rate $c(t)AB$ then the propensity is of this form.

In this common case (\ref{eq:special:form}) we have
\begin{equation}\begin{aligned}
	\theta_2 &= (\Delta t)\cdot A_c + \sum_{l=1}^r\left(g_l \cdot A_l\right),  \nonumber \\
	 \theta_4 &= (\Delta t)^2 \left(\sum_{l=1}^r \left(h_l  \cdot[A_l,A_c]\right) +
	 \sum_{l_1=1}^{r-1} \sum_{l_2=l_1+1}^r\left(g_{l_2}h_{l_1}-g_{l_1}h_{l_2}\right)[A_{l_1},A_{l_2}]\right).
\end{aligned}\end{equation}
where the scalar integrals are
\begin{equation}\begin{aligned}
	g_l =  \frac{1}{\Delta t}\int_{0}^{\Delta t}f_l(t+\tau) \du \tau,\quad h_l = \int_{-\Delta  t/2}^{\Delta t/2} \tau^j f_l(t+\Delta t/2+\tau) \du \tau, \quad l=1,\ldots,r.
\end{aligned}\end{equation}
Thus, the constant commutators $[A_c,A_l]$ and $[A_{l_1}, A_{l_2}]$ can be precomputed once and for all at the beginning.
Then, in each time step, only \textit{scalar} integrals are required.

\paragraph*{Magnus integrators are accurate for master equations, even for large time steps}
For (\ref{eq:magnus}) to be valid, Moan and Niesen have shown that the following inequality is a sufficient (but not necessary) condition \cite{Moan08},
\begin{equation}\label{eq:moanniesen}
	\int_{t}^{t+\Delta t}\| A(\tau)\|_{\ell_2} \du \tau \leq \pi.
\end{equation}
It is natural to wonder if the step size of a numerical method based on the Magnus expansion must be small enough to satisfy \eqref{eq:moanniesen}.
(That would place an unacceptably severe restriction on the step size whenever the matrix has a large norm, which is often the case for large sparse matrices arising in master equations.)
In practice, we find that we can solve the master equation with much larger time steps whilst maintaining stability (Figure~\ref{fig:cme20_moanniesen}).
This good experience is known for the time-dependent Schr{\"o}dinger equation and Hochbruck and Lubich \cite{Hochbruck03} provide an explanation. Unfortunately, the same reasoning does not explain this good experience in the context of the master equation.
 Their reasoning is based on the fact that the Hamiltonian of the Schr{\"o}dinger equation is self-adjoint and thus has only imaginary eigenvalues and can be transformed into a diagonal matrix with an orthogonal transformation.
 For the master equation, on the other hand, we usually do not have such a transformation.

\section{An estimate for state space truncated Magnus--Arnoldi propagation}\label{sec:full_estimate}

The numerical simulation with the Magnus--Arnoldi propagator based on a truncated state space has three sources of errors: the truncation of the state space, the truncation of the Magnus expansion, and the projection to the Krylov subspace.
We view the discretization process in three steps and use the following notation: First we truncate the state space and denote the solution of the corresponding perturbed master equation by $p^E$.
As a next step, we introduce a temporal grid, $0=t_0,t_1,\ldots,t_{N_t}=\tfin$ and apply the truncated Magnus expansion for propagation on each interval.
The solution of the corresponding perturbed master equation is denoted by $p^{EM}$.
Finally, we use the Arnoldi algorithm to compute the matrix exponential and get the fully perturbed master equation with solution $p^{EMA}$.

As usual in numerical methods for ODEs, we only have estimates for the local error over a small time step, so we ``glue'' these local estimates together over many time steps, to arrive at a global estimate at the end of the full time interval.
First we derive an estimate for a single time step.
The restarted perturbed problems then have an initial value that is additionally perturbed due to  errors from  previous time steps.

At each time step $n \in \{1,\ldots,N_t-1\}$, we have the  perturbed problem
\begin{equation}\begin{aligned}
	&\dot p^{EMA} = A p^{EMA} - A_R p^{EMA} + r_M(t), \quad t \in (t_n,t_n+1)\\
	&p^{EMA}(t_n) = p(t_n) + \rho_A(t_n)+e(t_n).
\end{aligned}\end{equation}
Here, $e(t_n) = p^{EMA}(t_n)-p(t_n)$ is the error accumulated from previous time steps and $-A_Rp^{EMA}$, $r_M$, and $\rho_A$ are perturbations due to state-space truncation, Magnus expansion truncation and the Arnoldi method, respectively.
Perturbations due to the Magnus expansion, $r_M$, can be estimated by (\ref{eq:magnus_residual}) and the perturbations due to the Arnoldi algorithm by (\ref{eq:arnoldi_residual}).
Note these residuals are computed with the truncated state space matrix $A_T$.

Applying the \textit{a posteriori} error estimation theory from Sec.~\ref{sec:apost} to the short time interval, we get
\begin{equation}
q(t_{n+1})^T e(t_{n+1}) = \int_{t_n}^{t_{n+1}} q(t)^T \left(r_M(t)-A_R p^{EMA}(t)\right) \du t + q(t_n)^T (\rho_A(t_n)+e(t_n)).
\end{equation}
We picked a special functional of the error that is based on the dual problem at time $t_{n+1}$.
In this way, the estimate from the previous step is the last term in the error estimate for the present step, $q(t_n)^T e(t_n)$, that covers errors from previous steps.
So we get a chain of estimates,
\begin{equation}
\begin{aligned}
	q(t_1)^Te(t_1) &= \int_{t_0}^{t_{1}} q(t)^T \left(r_M(t)-A_R p^{EMA}(t)\right)  \du t + q(t_0)^T (\rho_A(t_0)),\\
	q(t_{n+1})^T e(t_{n+1}) &= \int_{t_n}^{t_{n+1}} q(t)^T \left(r_M(t)-A_R p^{EMA}(t)\right)  \du t + q(t_n)^T (\rho_A(t_n)+e(t_n)),\\
	& n=1,\ldots,N_t-2,\\
	z^T e(t_{N_t}) &=  \int_{t_{N_t-1}}^{t_{N_t}} q(t)^T \left(r_M(t)-A_R p^{EMA}(t)\right)  \du t\\
	& + q(t_{N_t-1})^T (\rho_A(t_{N_t-1})+e(t_{N_t-1})).
\end{aligned}
\end{equation}

The local error in each step is
\begin{equation}
	 \int_{t_n}^{t_{n+1}} q(t)^T \left(r_M(t)-A_R p^{EMA}(t)\right)  \du t + q(t_n)^T \rho_A(t_n).
\end{equation}
Corresponding to (\ref{eq:est}), (\ref{eq:est1}), and (\ref{eq:est2}), we have the following three estimates:
\begin{equation}\label{eq:estdual}
	z^Te(t) \leq \sum_{n=0}^{N_t-1} \left( \int_{t_n}^{t_{n+1}} q(t)^T \left(r_M(t)-A_R p^{EMA}(t)\right)  \du t + q(t_n)^T \rho_A(t_n)\right),
\end{equation}

\begin{equation}\begin{aligned}\label{eq:estCSdual}
	|z^Te(t)| \leq&  \sum_{n=0}^{N_t-1} \left( \int_{t_n}^{t_{n+1}} \|q(t)\|_{\infty} \|r_M(t)\|_{\ell_1} \du t \right.\\ &\left.+\int_{t_n}^{t_{n+1}} \|q(t)\|_{\infty} \|A_R p^{EMA}(t)\|_{\ell_1} \du t + \|q(t_n)\|_{\infty}\| \rho_A(t_n)\|_{\ell_1}\right),
\end{aligned}\end{equation}
and
\begin{equation}\label{eq:estCS}\begin{aligned}
	& |z^Te(t)| \leq \\
	  &  \|z\|_{\ell_{\infty}} \sum_{n=0}^{N_t-1} \left( \int_{t_n}^{t_{n+1}}  \|r_M(t)\|_{\ell_1} \du t + \int_{t_n}^{t_{n+1}}  \|A_R p^{EMA}(t)\|_{\ell_1} \du t + \| \rho_A(t_n)\|_{\ell_1}\right).
\end{aligned}
\end{equation}

\section{An adaptive simulation algorithm}\label{sec:adaptive}

In the previous section, we derived an \textit{a posteriori} error estimate for propagation of the master equation.
We now explain how this estimate can be used to control the error by \textit{adaptively} choosing: the size of the state space, the size of the time step, and the dimension of the Krylov space.
We have other parameters to control the state space truncation error and the Arnoldi error.
Therefore, we primarily choose the step size to control the Magnus error.
Then, we choose the state space as large as needed for the chosen step size, and choose the Krylov subspace large enough to keep the Arnoldi error small.
This strategy works well when the time-dependent nature of the problem requires small time steps.
In general, Krylov methods usually only give good approximations for a moderate dimension (say, 40) of the Krylov space.
Therefore, we also aim to keep the time step small enough to keep the size of the Krylov space  moderate.
In a first approximation,  $\rho_A$ is proportional to $\Delta t$, so we may assume the ratio between $\Delta t$ and the size $s$ of the Krylov subspace is constant.
Given $s_n$ and $\Delta t_n$ for step $n$, we have the following upper bound for the size of $\Delta t_{n+1}$,
\begin{equation}
	\Delta t_{n+1,\max} = \frac{s_{\max}}{s_n} \Delta t_n,
\end{equation}
where $s_{\max}$ is the maximum acceptable dimension of the Krylov subspace.
Given a tolerance $\tol$, the simulation based on estimate (\ref{eq:estCS})  is performed as:
\begin{enumerate}
	\item \textit{Try to reduce state space}: Every $C$ steps (for some constant $C$) check if the state space can be reduced.
	\item \textit{Check if the state space needs to be increased}: Check whether or not $A_R(t_n+\frac{\Delta t_n}{2})\left(p(t_n)+p(t_{n+1})\right)/2< \frac{\tol}{\tfin}$. If not, increase the state space.
	\item \textit{Propagate solution with adaptive Arnoldi method}: Propagate $n$th step.
	Choose size of Krylov subspace such that $\rho_{A}(t_n) < \frac{\tol\cdot \Delta t_n}{\tfin}$.
	\item \textit{Compute new step size according to Magnus residual}: Estimate the error due to the truncation of the Magnus expansion via \eqref{eq:magnus_residual} and compute the optimal time step for the next step according to the formula,
	\begin{equation}
		\Delta t_{n+1}^M = \left( \Delta t_n \cdot \frac{\tol / \tfin}{\| \theta_{m+2}(t_{n-1}+\Delta t_n)\cdot p(t_n)\|_{\ell_1}}\right)^{1/m} \cdot \Delta t_n.
	\end{equation}
	\item \textit{Keep the step size small enough for the Arnoldi method}: Check whether or not the step size $\Delta t_{n+1}^M$ is acceptable and set
	\begin{equation}
		\Delta t_{n+1} = \min \left(\Delta t_{n+1}^M,\frac{s_{\max}}{s_n} \Delta t_n\right).
	\end{equation}
\end{enumerate}

If we want to use estimate (\ref{eq:estdual}) or (\ref{eq:estCSdual}) instead of (\ref{eq:estCS}), then we additionally have to include the dual or its $\ell_{\infty}$ norm, respectively, in the error checking.
Note that Step 4 is redundant if the  matrix is constant, in which case we set $\Delta t_{n+1} = \frac{s_{\max}}{s_n} \Delta t_n$.



\section{Dual problem}\label{sec:dual}

Estimates (\ref{eq:estdual}) and (\ref{eq:estCSdual}) require the solution of a dual problem while estimate (\ref{eq:estCS}) only requires the primal solution.
The more information we use from the dual problem, the more accurate our estimate of the error will be.
Thus larger --- and fewer --- time steps can be used.
On the other hand, solving the dual problem is a computational cost.
Hence, there is a trade-off between costs and benefits of computing the dual.
Such trade-offs are familiar in the literature on adjoint methods, arising much more generally than merely in our applications.

Memory is an additional difficulty.
 For the estimates (\ref{eq:estdual}) and (\ref{eq:estCSdual}) we need primal and dual solutions at the same time.
 Therefore, we have to first solve the dual problem backwards in time, and then save this dual solution in memory, to then be ready to compute the solution of the forward problem.
Note estimate (\ref{eq:estCSdual}) requires only the scalar $\|q(t)\|_{\infty}$ (and not the vector $q(t)$), which will usually fit in memory.
Estimate (\ref{eq:estdual}) needs the full solution so it requires much more memory or checkpointing with additional requirements on computation.

To be competitive with estimate (\ref{eq:estCS}), we cannot spend too much effort on the dual problem so we only compute a rough approximation.
However, we will, of course, not know a priori how the error in the dual behaves.
In principle, we can appeal to the same theory as for the primal.
Since we do not have the solution of the `dual of the dual' --- i.e., the primal --- available, we  may resort to using an estimate corresponding to (\ref{eq:estCS}) for the dual.
Denoting the error in the dual solution at time zero by $e_q$ and the residual by $r_q$, we find the following estimate for a functional defined by $z$,
\begin{equation}\begin{aligned}
	|z^Te_q| &\leq \int_{\tfin}^0 p(t)^Tr_q(t) \du t + p(\tfin)^Tq(\tfin) \\
	&\leq \int_{\tfin}^0 \|p(t)\|_{\ell_1}\|r_q(t)\|_{\ell_{\infty}} \du t +\|p(\tfin)\|_{\ell_1}\|r_q(\tfin)\|_{\ell_{\infty}}
\\&	 \leq \int_{\tfin}^0 \|r_q(t)\|_{\ell_{\infty}} \du t+\|r_q(\tfin)\|_{\ell_{\infty}}.
\end{aligned}\end{equation}
So if we make sure that
\begin{equation}
	\max_{t_n \leq t\leq t_{n+1}}\|r_q(t)\|_{\infty} \leq \frac{\varepsilon}{\tfin} \Delta t
\end{equation}
in each time step, then we bound the error at any time by $\varepsilon$.
Of course, we should choose  $\varepsilon$ to be much larger than the tolerance to which we want to compute the solution.
Choosing $\varepsilon$ is therefore a familiar trade-off between the advantage of the more accurate error estimates that come from more accurate dual solutions, and the disadvantage of the cost of solving the dual.


Altogether, we conclude that estimate (\ref{eq:estCS}) is the most reliable estimate, because it does not depend on uncontrolled influences from the approximation of the dual.
We invest extra effort solving the primal problem more accurately (in the subspace perpendicular to $z$) instead of solving an additional dual problem.
This is the strategy of choice when $z$ is unknown and also often when we are interested in several functionals at the same time.

\section{Numerical experiments}\label{sec:numerical}

We now present numerical results to complement our theoretical findings.
Given a prescribed tolerance for the error in a given component, we simulate the solution of the master equation with our state-space truncated Magnus--Arnoldi propagator with adaptivity based on the three estimates.
We compare effectiveness and efficiency of the three estimates.
By effectiveness, we mean that the error is indeed below the tolerance.
The closer the actual error is to the tolerance, the more efficient the estimate.

\subsection{A test case with analytic solution}

Consider a master equation, with time varying propensities,
\begin{equation}\label{eq:2d}\begin{aligned}
	&\frac{\du}{\du t} p = A(t)p = \left[\left(\begin{array}{rr}  -1& 1\\1& -1
	\end{array}\right)+ \sin(t)\left(\begin{array}{rr}  -1& -1\\1& 1
	\end{array}\right) \right]
	 p, \\
	 &p(0) = (\sigma,1-\sigma)^T,
\end{aligned}\end{equation}
for $\sigma \in [0,1]$.
The solution is
\begin{equation}\begin{aligned}\label{eq:example1_analytic}
	p(t) = \left( \begin{array}{c}
\frac{1}{2}+\frac{1}{5}\cos(t)-\frac{2}{5}\sin(t)+\left(\sigma-\frac{7}{10}\right)\exp(-2t)\\
	\frac{1}{2}-\frac{1}{5}\cos(t)+\frac{2}{5}\sin(t)-\left(\sigma-\frac{7}{10}\right)\exp(-2t)
	\end{array}\right).
\end{aligned}\end{equation}

Table~\ref{tab:2d_10} shows results for simulations of equation (\ref{eq:2d}) with $\sigma=1$ over a time interval of $t \in [0,10]$.
The experiments show that it does not pay off to compute the dual solution for this example.
Table~\ref{tab:2d_10} shows all estimates are effective.
It also shows the estimates are more efficient as more information from the dual is used, as one would expect from theory.
Note that in this example, the only source of error is due to truncation in the Magnus expansion.

\begin{table}[htdp]
\caption{Results for the simple master equation in (\ref{eq:2d}) with the error checked by comparison of the numerical solution with the exact, analytic solution. Parameters: $t_{\textrm{final}} =10$ and tolerance $ = 10^{-3}$.}
\begin{center}
\begin{tabular}{|c|c|c|}
\hline
estimate & $\ell_{\infty}$ error & \# time steps (+dual steps) \\
\hline
(\ref{eq:estdual}) &$9.1\times 10^{-4}$  & 62+100\\
(\ref{eq:estCSdual})&$3.6\times 10^{-4}$ &112+100\\
(\ref{eq:estCS})&$3.3\times 10^{-4}$ & 131\\
\hline
\end{tabular}
\end{center}
\label{tab:2d_10}
\end{table}%

\subsection{Chemical master equation}\label{sec:CME}


Reversible isomerization involves the conversion of a molecule between two different forms.
Consider a chemical master equation (CME) describing this process:
\begin{equation} \label{eq:isomerization}
  X \xrightarrow{c_{12}} Y, \qquad Y \xrightarrow{c_{21}} X.
\end{equation}

We choose $c_{1,2}(t) = 1+\sin(t)$ and $c_{2,1}(t) = 1-\sin(t)$ and let the initial distribution be a binomial distribution $\mathcal{B}(x_k;N,p_0)$ with $N=2000$ molecules and parameter $p_0 = (1/3,\, 2/3)^T$.
Observing the conservation law, there are $2001$ possible states $x_k = (k,N-k)$ for $k=0,\ldots,2000$.
The chemical master equation for \eqref{eq:isomerization} is,
\begin{equation}\begin{aligned}\label{eq:CME}
	\frac{\du}{\du t} P(t) = \left(A_0+ \sin(t) A_1\right) P(t),
\end{aligned}\end{equation}
where the $k$th element of $P$ gives the probability that the system is in state $x_k$ and
\begin{equation}\begin{aligned}
	A_0^{(i,j)} = \left\{ \begin{array}{l  l}
	-N & \textrm{if } i=j\\
	 j & \textrm{if } i=j-1\\
	 N-j & \textrm{if } i=j+1\\
	 0 & \textrm{otherwise}\\
	\end{array}\right. , \text{ and }
	A_1^{(i,j)} = \left\{ \begin{array}{l  l}
	N-2j & \textrm{if } i=j\\
	 j & \textrm{if } i=j-1\\
	 -N+j & \textrm{if } i=j+1\\
	 0 & \textrm{otherwise}\\
	\end{array}\right. .
\end{aligned}\end{equation}
The exact analytic solution to (\ref{eq:CME}) is known (see Proposition 1 in \cite{Jahnke07}) to be
\begin{equation}
	P_k(t) = \mathcal{B}(x_k;N,p(t)), \quad k=0,\ldots,N,
\end{equation}
where the $1 \times 2$ parameter vector $p(t)$ is the solution of the reaction rate equation
\begin{equation}\label{eq:reaction_rate_equation}
	\dot p(t) = A(t)p(t), \qquad p(0) = p_0,
\end{equation}
in which $A(t)$ is as in (\ref{eq:2d}) of the previous example with $\sigma=1/3$ and thus expression (\ref{eq:example1_analytic}) gives the solution $p(t)$.

Table~\ref{tab:CME} compares the control strategies based on estimates (\ref{eq:estdual}),  (\ref{eq:estCSdual}), and (\ref{eq:estCS}).
As the error functional, we use the error in component 1150, where the solution is relatively large but not maximal.
In this example, the time dependence is weak.
The step size control tends to lead to Krylov subspaces of successively larger dimension.
In this case, we need a step size control that also takes into account the Arnoldi residual.
So step 5 in our algorithm becomes relevant.

The tolerance $10^{-5}$ might seem at first to be too small, but bear in mind that it concerns the absolute error, and the value of the solution in component 1050 is $1.5 \times 10^{-3}$. The dual solution has been computed with step size control based on an estimate similar to (\ref{eq:est2}) for the dual master equation with a tolerance of $10^{-4}$.
We see that it does not pay off to compute the dual in this example.
 Looking at the error in component 1150, we see that all estimates are effective and the efficiency increases as more information from the dual is included.
  As expected from theory, estimate (\ref{eq:estCS}) also ensures that the $\ell_{\infty}$ error is below the tolerance.

\begin{figure}
\begin{center}
\includegraphics[scale=0.9]{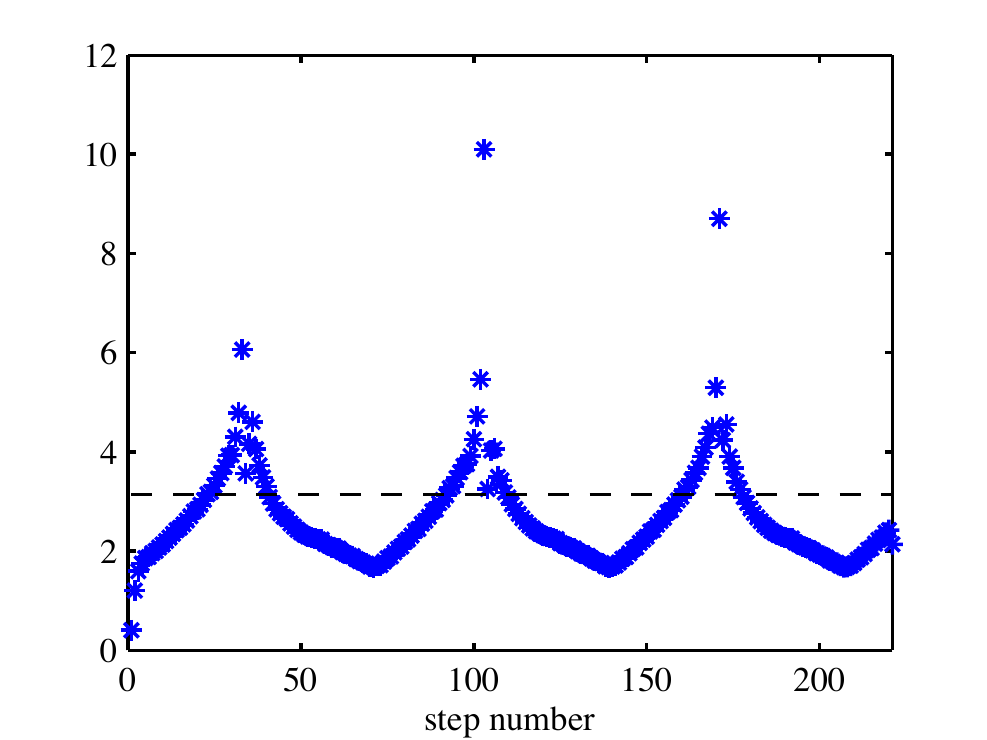}
\caption{Value of $\left\|\int_{t}^{t+\Delta t} A(\tau) \du \tau\right\|_{\ell_2}$ for each time step (***), upper bound for convergence in Magnus expansion due to Moan and Niesen (-- -- --).
CME example from Sec.~\ref{sec:CME}.}
\label{fig:cme20_moanniesen}
\end{center}
\end{figure}

In this experiment, the Moan--Niesen criterion \eqref{eq:moanniesen} is violated in each time step.
Nevertheless, we get an accurate result.
This is numerical evidence  that large time steps with a Magnus based integrator for the explicitly time-dependent master equation can still be very accurate.
However, in this example, with $N = 2000$, the error due to the Magnus-expansion is almost insignificant compared to the Arnoldi error.
Therefore, we repeated the experiment with $N = 20$.
In this case, the Magnus error controls the step size and not the Arnoldi error.
For each time step, we compute $\left\|\int_{t}^{t+\Delta t} A(\tau) \du \tau\right\|_{\ell_2}$ because this is a lower bound for $\int_{t}^{t+\Delta t}\| A(\tau)\|_{\ell_2} \du \tau$.
Fig.~\ref{fig:cme20_moanniesen} shows an estimate of the Moan-Niesen criterion for a simulation with tolerance $10^{-3}$.
Again, in this experiment the solution is recovered correctly ($\ell_{\infty}$ error of the size $1.3 \times 10^{-4}$) even though the Moan-Niesen criterion is violated in several steps.

\begin{table}[htdp]
\caption{Errors and necessary number of matrix-vector products (MVPs) for $t_{\textrm{final}} =10$ and tolerance $10^{-5}$ for the CME with varying propensity (Sec.~\ref{sec:CME}).}
\begin{center}
\begin{tabular}{|c|c|c|c|}
\hline
estimate & $\ell_{\infty}$ error & error in  comp.~1150 &\# MVPs (primal+dual) \\
\hline
(\ref{eq:estdual})&$1.5\cdot10^{-4}$ & $6.4\cdot10^{-6}$ &28201+27191\\
(\ref{eq:estCSdual})&$1.0\cdot10^{-5}$ & $3.1\cdot10^{-6}$ &21800+27191\\
(\ref{eq:estCS})&$8.1\cdot10^{-7}$ & $2.5\cdot10^{-7}$ &31928\\
\hline
\end{tabular}
\end{center}
\label{tab:CME}
\end{table}%


\begin{figure}
\begin{center}
\includegraphics[scale=0.9]{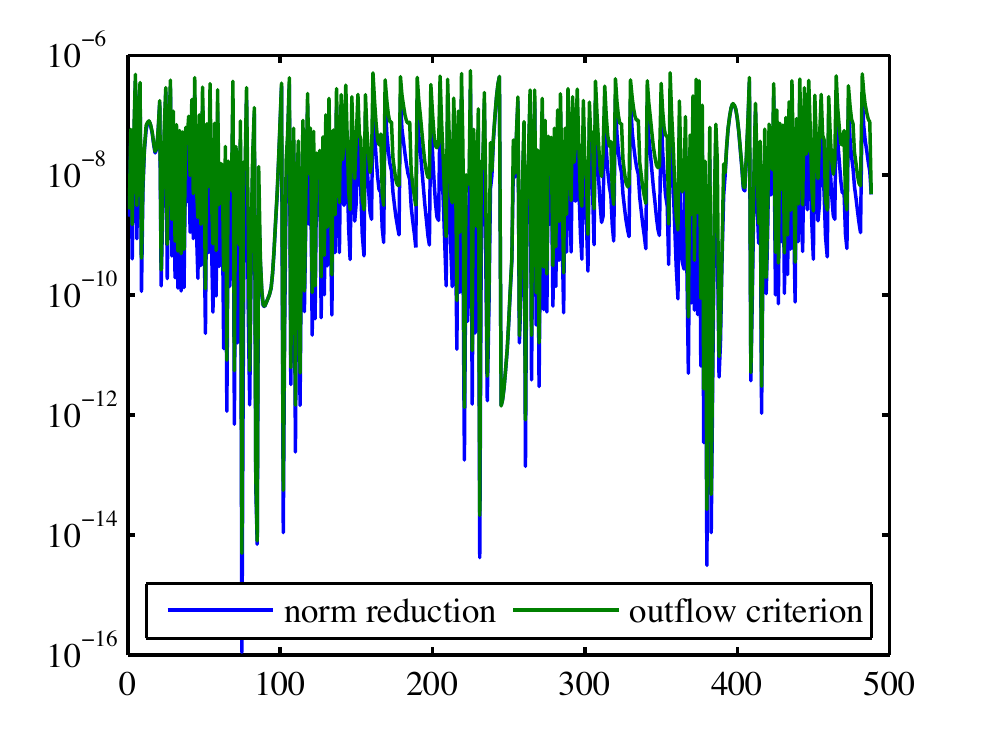}
\caption{Comparison of our (green) outflow-criterion  \eqref{eq:res_sst} with the (blue) criterion \eqref{eq:FSP:est} of Munsky and Khammash \cite{Munsky06} based on the loss of probability for estimating the error due to state space truncation.
CME example from Sec.~\ref{sec:CME} (without a Krylov approximation).}
\label{fig:CME_SST}
\end{center}
\end{figure}

So far, we have computed on the whole state space.
However, most of the probability mass of the solution is confined to a certain region, which is oscillating in time.
We therefore repeated the experiment, but this time we allowed truncation of the state space.
We start with 400 states.
The adaptive control of the state space successively reduces this to about 250 states.
After this initial phase, the size varies between 121 and 261 states.
In this case, the problem is small enough that we can use the Matlab function \texttt{expm} to compute the matrix exponential.
This (\texttt{expm}) almost conserves the $\ell_1$ norm of the solution.
Hence, we could just as well use the `probability-loss' criterion  \eqref{eq:FSP:est} proposed by \cite{Munsky06} to control the error due to state space truncation (cf.~Sec.~\ref{sec:statespace}).
Figure~\ref{fig:CME_SST} compares both criteria  \eqref{eq:res_sst} and \eqref{eq:FSP:est}.
 Both criteria show  very similar behavior but the outflow criterion \eqref{eq:res_sst} is slightly overestimating in most steps.
 This shows that our criterion \eqref{eq:res_sst}  is comparable to an estimate from the literature in examples where the latter can be applied.

\subsection{CME with constant propensities}\label{sec:cme_const}

Our adaptive algorithm also applies to examples with a constant matrix $A$.
Table~\ref{tab:CME_const} summarizes results for the same example \eqref{eq:CME}, except that we set $A=A_0$.
The algorithm behaves similarly.

\begin{table}[htdp]
\caption{Errors and necessary number of matrix-vector products (MVPs) for $t_{\textrm{final}} =10$ and tolerance $10^{-5}$ for the CME with constant propensity matrix (Sec.~\ref{sec:cme_const}).}
\begin{center}
\begin{tabular}{|c|c|c|c|}
\hline
estimate & $\ell_{\infty}$ error & error in  comp.~1050 &\# MVPs (primal+dual) \\
\hline
(\ref{eq:estdual})&$7.6\cdot10^{-5}$ & $2.4\cdot10^{-5}$ &1607+11364\\
(\ref{eq:estCSdual})&$2.2\cdot10^{-6}$ & $2.1\cdot10^{-7}$ & 1941+11364\\
(\ref{eq:estCS})&$4.6\cdot10^{-10}$ & $7.0\cdot10^{-11}$ &2366\\
\hline
\end{tabular}
\end{center}
\label{tab:CME_const}
\end{table}%

\subsection{T Cell Homeostasis}\label{sec:tcell}

Finally, we demonstrate our approach on a larger example, modeling T cell homeostasis \cite{Stirk08,MacNamara10}.
The master equation has time-dependent propensities
\begin{equation}
\begin{aligned}
\frac{\ud}{\ud t} p(n,n',t) &= \lambda_{n-1,n'} p(n-1,n',t)+\lambda'_{n,n'-1} p(n,n'-1,t)\\
	&+\mu_{n+1,n'} p(n+1,n',t)+\mu'_{n,n'+1} p(n,n'+1,t)\\
	&-(\lambda_{n,n'}+\lambda'_{n,n'}+\mu_{n,n'}+\mu'_{n,n'}) p(n,n',t),\\
p(n,n',0) &= \delta_{n 10} \delta_{10 n'}.
\end{aligned}
\end{equation}
Here, $\mu_{n,n'} = n$, $\mu_{n,n'}' = n'$, $\lambda_{n,n'} = \lambda_{n,n'} (t)= 30 n \left(\frac{1}{n+n'}+ \frac{1}{n+1000} \right) \frac{1}{1+\left(\frac{t}{15×}\right)^5}$ and $\lambda_{n,n'}' = \lambda_{n,n'}' (t)= 30 n \left(\frac{1}{n+n'}+ \frac{1}{n'+1000} \right) \frac{1}{1+\left(\frac{t}{15×}\right)^5}$.

\begin{figure}
  \centering
  \subfigure[Step size (red) and Krylov space (black)]{
    \label{subfig:immunology_krylov_deltat}
    \includegraphics[width=\textwidth]{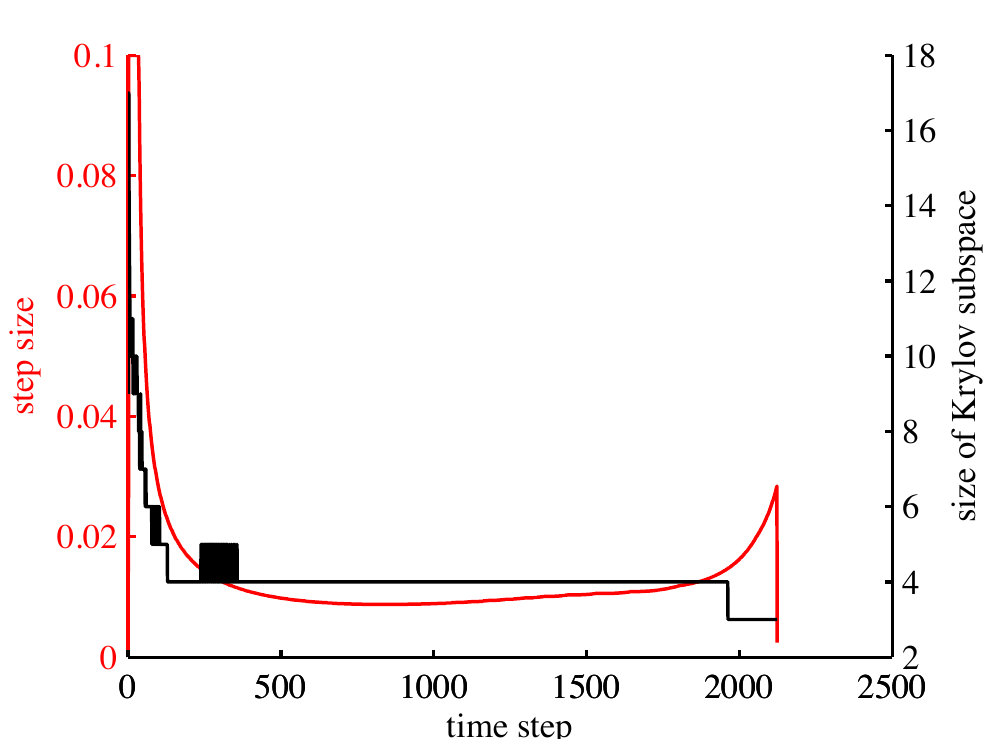}}
  \subfigure[Size of state space]{
    \label{subfig:immunology_statespace}
    \includegraphics[width=\textwidth]{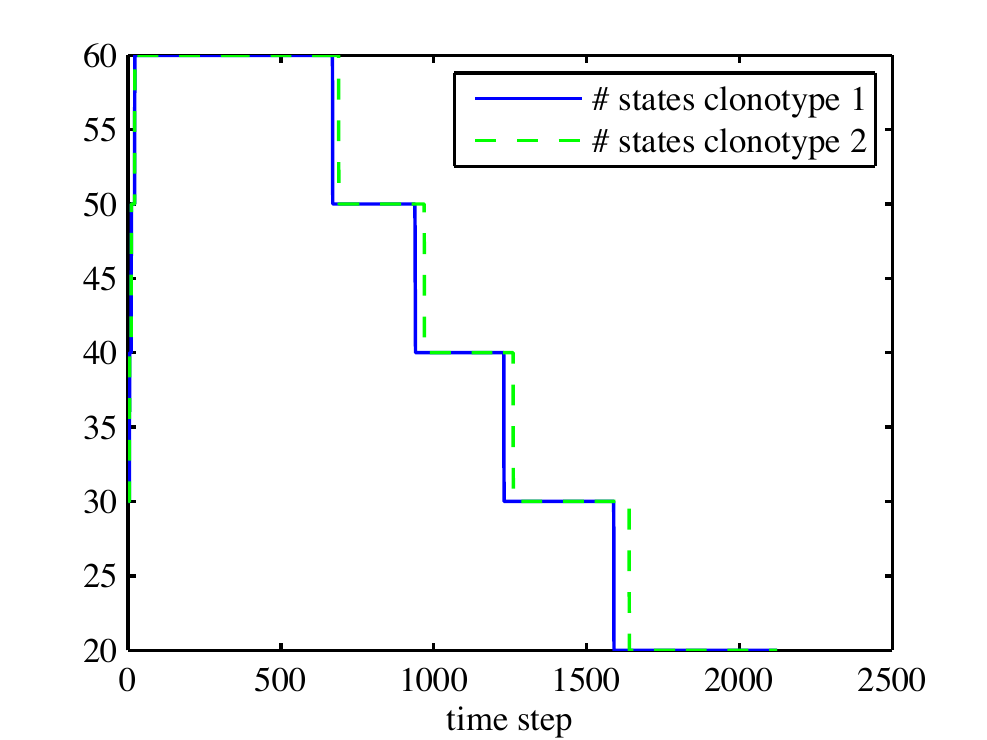}}
  \caption{Choice of temporal step size, dimensionality of the Krylov subspace, and  the state space for the immunology example in Sec.~\ref{sec:tcell}.}
\label{fig:immunology}
\end{figure}

For this example, we also truncate the state space.
We start with a state space that includes 0 to 29 copies of each species.
The size grows quickly in the beginning, to a maximum of 59 copies of each species.
Then, the size of the state space is gradually reduced with time.
Figure~\ref{subfig:immunology_statespace} shows the size of the state space as a function of the time step.
Figure~\ref{subfig:immunology_krylov_deltat} shows how the size of the time step changes.
It is closely related to the nature of the time-dependence in the propensity matrix, which is reasonable given the expression for the error term due to the Magnus expansion.
We do observe some kind of drift  towards smaller time steps.
This is probably due to starting from a simple Dirac delta distribution.
As time evolves, the shape of the probability distribution gets more complicated, which is why the Krylov algorithm will have to do more work.
A similar behaviour is also seen in the other  examples above.
The number of Arnoldi iterations in each time step is nearly constant  (Figure~\ref{subfig:immunology_krylov_deltat}).
Near the end of the simulation, the Krylov subspace is very small even though the time steps are getting a little larger.
This shows that the size of the Krylov space is not only related to the size of the time step, but also to the nature of the time-dependence.

We do not have an exact analytic solution for this example, so we use a computation with a large state space and small time steps as a `reference solution.'
In this way, we estimate the $\ell_{\infty}$ error of our simulation to be $2.3 \cdot 10^{-5}$.
This has been achieved with 2124 time steps and 8704 matrix-vector products.
If instead we take equally many time steps with a fixed step size, but allow (adaptively) variable sizes of both the state space and of the Krylov subspace, then we get an $\ell_{\infty}$ error of $7.9 \cdot 10^{-5}$, with 8289 matrix-vector products.
Fixing also the size of the Krylov subspace, yields an error of $3.4 \cdot 10^{-3}$ for 8496 matrix-vector products (four per time step), or $3.5 \cdot 10^{-4}$ for 10620 matrix-vector products (five per time step).
A benefit of the framework for \textit{a posteriori} error estimation is that it makes adaptivity possible.
These experiments demonstrate that allowing adaptivity in both the step size and in the Krylov subspace yields a considerably reduced error for the same amount of work.



%




\section{Discussion}\label{sec:discussion}
We have devised and compared three types of \textit{a posteriori} error estimates for the master equation.
Numerical tests show that the most pessimistic estimate, \eqref{eq:est2}, which does not require a numerical solution of the dual problem, is the most robust.
It is also the most general because it provides an error estimate for many functionals of the error at the same time, as in \eqref{eq:simple:error:bound}.
Often, it is the most computationally efficient.

When choosing a numerical method it is desirable to mimic the conservation properties of the continuous problem.
With a matrix $A$ of the type that arises in master equations, we have $\ell_1$ norm conservation, as well as positivity of the vector $p$.
That is, $e^Tp(t) \equiv 1$ and $p(t) \geq 0$ for all $t$.
The first property is a linear invariant and is therefore preserved by any Runge--Kutta method \cite{Gear92,HaiLub06}.
However, most Runge--Kutta methods do not respect positivity.
Indeed, it has been shown that unconditionally positive Runge--Kutta methods can only be first order \cite{Bolley78}.
Uniformization is an alternative but also only of low order \cite{SidjeUniformization07}.
There have been efforts to devise other second-order schemes \cite{Bruggeman07, Broekhuizen08}.
We found this type of scheme inefficient for our examples.

The Magnus--Arnoldi propagator does not fully preserve the geometric stuctures of the master equation.
These structures are coupled to the two properties \eqref{eq:zero:column:sum} and \eqref{eq:nonnegative:off:diagonal} of a matrix $A(t)$ for the master equation.
It would be nice if the truncated Magnus propagator shared these same properties.
It is easy to verify that the first property of `zero column sum' is conserved regardless the number of terms included in the approximation of the Magnus propagator.
However, commutation of matrices does not preserve the second property, which is why truncated Magnus propagators of order greater than two may have negative off-diagonal entries (unlike $A(t)$).
Finally, the Arnoldi process conserves the $\ell_2$ but not the $\ell_1$ norm of a vector.
In practice, we observe that the norm-conservation error is of considerably smaller size than the $\ell_2$- or $\ell_{\infty}$-error of the solution.

Alternatively, one can use a numerical integrator that does not mimic the conservation properties of the continuous problem and simply project the propagated solution onto the manifold defined by the conservation properties (cf.~\cite{Sandu01}).
However, this strategy is not suitable when a truncated space is used in combination with estimates such as
\eqref{eq:FSP:est} of \cite{Munsky06} that exploit a loss of probability to estimate the error due to state space truncation.

\section{Conclusions}\label{sec:conclusions}
An \textit{a posteriori} error estimate for the master equation, based on adjoint methods, has been derived.
As a general purpose estimate, it can be combined with many different types of numerical approaches to the master equation.
As one example, we used the estimate to develop an adaptive algorithm that incorporates a Magnus--Arnoldi propagator into this framework.
This type of evolution operator has proven to be efficient for the master equations in immunology and in chemistry that we have studied.
Part of the reason for this is that relatively large time steps can be used whilst maintaining good accuracy and stability.
Previous methods have employed error estimates, associated with truncation of the state space, that are valid when it is assumed that the matrix exponential is computed exactly.
This is an issue because in practice the exponential is almost always approximated.
The adjoint framework has allowed us to address this issue by suggesting an alternative   error estimate that can be incorporated into the framework in a way that is compatible with a Krylov approximation to the matrix exponential.
Overall, adjoint methods offer a mathematical framework for combining and controlling errors from different sources, giving more confidence in the final computed solution.

\section{Acknowledgment}
Katharina Kormann began this work in Uppsala University.

\bibliography{master}

\end{document}